\documentclass[12pt]{article}

\usepackage{amssymb}

\newtheorem {theorem} {Theorem}
\newtheorem {definition} [theorem]{Definition}

\newcommand{\bbox}{\ \hfill\rule[-1mm]{2mm}{3.2mm}}

\def\div{{\rm div}}

\title{A note on: ``Relaxation Oscillators with Exact Limit Cycles''.\thanks{The authors are
partially supported by a MCYT grant number BFM 2002-04236-C02-01.
The first author is also partially supported by DURSI of
Government of Catalonia ``Distinci\'o de la Generalitat de
Catalunya per a la promoci\'o de la recerca universit\`aria''.}}

\author{{\sc Jaume Gin\'e and Maite Grau}}
\date{}

\begin{document}
\maketitle \vspace{-1.0cm}

\begin{center}
Departament de Matem\`atica. Universitat de Lleida. \\
Avda. Jaume II, 69. 25001 Lleida, SPAIN. \\ {\rm E--mails:} {\tt
gine@eps.udl.es}, {\tt mtgrau@matematica.udl.es}
\end{center}
\begin{abstract}
In this note we give a family of planar polynomial differential
systems with a prescribed hyperbolic limit cycle. This family
constitutes a corrected and wider version of an example given in
the work \cite{Abdelkader}. The result given in this note may be
used to construct models of Li\'enard differential equations
exhibiting a desired limit cycle.
\end{abstract}

{\small{\noindent 2000 {\it AMS Subject Classification:} 34A05, 34C07.   \\
\noindent {\it Key words and phrases:} Planar polynomial
differential system, limit cycle, hyperbolicity, Li\'enard
differential equation.}}

\section{Introduction and statement of the main result \label{sect1}}

Our purpose in this work is to give a family of planar polynomial
differential systems of the form:
\begin{equation}
\dot{x}= P(x,y), \quad  \dot{y}=Q(x,y), \label{eq1}
\end{equation}
for which an explicit expression of a limit cycle, that is, an
isolated periodic orbit, can be given. We assume that $P(x,y)$ and
$Q(x,y)$ belong to the ring of real polynomials in two variables
$\mathbb{R}[x,y]$, and we will always assume that $P(x,y)$ and
$Q(x,y)$ are coprime polynomials. We denote by ${\rm d}= \max \{
\deg P, \deg Q \}$ and we say that ${\rm d}$ is the degree of
system (\ref{eq1}).
\newline

In the work \cite{Abdelkader} a family of planar polynomial
differential systems like (\ref{eq1}) is studied and the existence
of an explicit limit cycle is pretended to be given. The author of
\cite{Abdelkader} gives a family of systems of the form
(\ref{eq1}) with a prescribed invariant algebraic curve. This
curve $f(x,y)=0$ has an oval surrounding the origin of
coordinates. However, in \cite{Abdelkader} there is no proof of
the fact that the oval of $f(x,y)=0$ is an isolated periodic
orbit, that is, a limit cycle. It is stated as obvious. We have
been able to weaken the hypothesis appearing in \cite{Abdelkader},
getting a bigger family of planar polynomial systems, and we have
been able to show that the oval of $f(x,y)=0$ is a hyperbolic
limit cycle.

\begin{theorem}
We consider a polynomial $p(x)$ such that $p(x_e)=0$ and
$p(x_d)=0$ for some $x_e<0$ and $x_d>0$, $p'(x_e) \neq 0$ and
$p'(x_d) \neq 0$. We assume that $p(x)>0$ for all $x$ in the
interval $(x_e,x_d)$. We consider another polynomial $q(x)$
satisfying $p(x) q(x)^2 \neq 1$ for all $x \in (x_e,x_d)$ and
$q'(x) \neq 0$ for all $x \in (x_e,x_d)$. Then, the algebraic
curve given by $f(x,y)=0$ with $f(x,y):= (y-p(x)q(x))^2 - p(x)$
has an oval in the band $x_e \leq x \leq x_d$ which is a
hyperbolic limit cycle for the following system:
\begin{equation}
\begin{array}{lll}
\displaystyle \dot{x} & = & \displaystyle y, \vspace{0.2cm} \\
\displaystyle \dot{y} & = & \displaystyle \left(\frac{3}{2} \,
q(x) p'(x) + p(x) q'(x) \right) y - \displaystyle \frac{p'(x)}{2}
\left( p(x) q(x)^2-1 \right).
\end{array} \label{eq2}
\end{equation} \label{th1}
\end{theorem}

\vspace{-0.3cm}

We note that a system of the form (\ref{eq2}) can also be viewed
as an autonomous Li\'enard differential equation:
\begin{equation}
\ddot{x} + f(x) \dot{x} + g(x) \,=0, \label{eq3}
\end{equation}
where $f(x)$ and $g(x)$ are the polynomials given by: $f(x):= - 3
q(x) p'(x)/2 - p(x) q'(x)$ and $g(x):= (p'(x)/2) \left( p(x)
q(x)^2-1 \right)$. Therefore, Theorem \ref{th1} may be used to
construct models of Li\'enard differential equations exhibiting a
hyperbolic limit cycle.
\par
One of the most famous Li\'enard differential equation is called
van der Pol equations and it appears when studying the vacuum-tube
circuits. This particular equation (\ref{eq3}) has $f(x)=\mu
(x^2-1)$ and $g(x)=1$, with $\mu \in \mathbb{R}$, and it exhibits
a unique hyperbolic limit cycle surrounding the origin. This limit
cycle is shown to be non-algebraic in the work of Odani
\cite{Odani}. We are not considering van der Pol's equation since
the systems described in Theorem \ref{th1} always exhibit an
algebraic limit cycle. The systems given in (\ref{eq2}) are
examples of Li\'enard equations with an algebraic and hyperbolic
limit cycle.
\par
The same kind equations are studied in the work \cite{Abdelkader},
but under other hypothesis for the polynomials $p(x)$ and $q(x)$,
and the author of \cite{Abdelkader} pretends to state the
existence of a limit cycle. The conditions for the polynomials
$p(x)$ and $q(x)$ appearing in \cite{Abdelkader} are: $p(x)$ is an
even polynomial, $p(0)>0$, there exists a value $X>0$ such that
$p(X)=0$, $p(x)<0$ for all $x>X$, $q(x)$ is an odd polynomial and
$q(0)=0$. All these conditions are contained in the ones that we
assume in Theorem \ref{th1}. However, the authors noticed that in
the work \cite{Abdelkader} the condition $p(x) q(x)^2 \neq 1$ for
$x \in (-X,X)$ does not appear and it is not implied by the other
hypothesis. As we will see in the proof of Theorem \ref{th1},
condition $p(x) q(x)^2 \neq 1$ for $x \in (x_e,x_d)$ is necessary
to have a limit cycle and it cannot be avoided.
\newline

We remark that each system of the family (\ref{eq2}) has an
algebraic limit cycle, which is an oval of the real algebraic
curve $f(x,y)=0$, and it may have other limit cycles which are not
taken into consideration. These other limit cycles can be
contained in $f(x,y)=0$ or not. If they are contained in an
invariant curve, we can treat them with the same methods described
in this note. For instance, it can be shown that the system
(\ref{eq2}) with $p(x)=(1-x^2)(4-x^2)(9-x^2)$ and $q(x)=x/100$ has
$3$ hyperbolic limit cycles all of them contained in the
corresponding invariant algebraic curve $f(x,y)=0$.
\newline

In order to prove Theorem \ref{th1}, we need some results relating
ovals of curves of the form $f(x,y)=0$ with the fact that the oval
of such a curve is a limit cycle. These preliminary results are
stated in Section \ref{sect2}. Once we have stated these previous
results, we prove Theorem \ref{th1} in Section \ref{sect3}.

\section{Preliminary results \label{sect2}}

We are considering limit cycles which are contained in a real
curve $f(x,y)=0$, which does not need to be algebraic. This fact
leads us to the definition of invariant of a system (\ref{eq1}).
\begin{definition}
Let us consider an open set $\mathcal{U} \subseteq \mathbb{R}^2$
and a $\mathcal{C}^1(\mathcal{U})$ real function denoted by
$f(x,y): \mathcal{U} \subseteq \mathbb{R}^2 \to \mathbb{R}$. We
say that $f(x,y)$ is an {\em invariant} for a system {\rm
(\ref{eq1})} if \begin{equation} P(x,y) \, \frac{\partial
f}{\partial x}(x,y) \, + \, Q(x,y) \, \frac{\partial f}{\partial
y} (x,y) \, = \, k(x,y) \, f(x,y), \label{eqci} \end{equation}
with $k(x,y)$ a polynomial of degree lower or equal than ${\rm
d}-1$, where {\rm d} is the degree of the system. This polynomial
$k(x,y)$ is called the {\em cofactor} of $f(x,y)$. \label{def3}
\end{definition}
In case that $f(x,y)$ is a polynomial we say that $f(x,y)=0$ is an
{\em invariant algebraic curve} for system (\ref{eq1}). We notice
that if $f(x,y)$ is an invariant of system (\ref{eq1}) and
$f(x,y)=0$ defines a curve in the real plane, then the function
$P(x,y) \, (\partial f/\partial x) \, + \, Q(x,y) \, (\partial
f/\partial y)$ equals zero on the points such that $f(x,y)=0$.
This fact implies that the real curve $f(x,y)=0$ is formed by
orbits of system (\ref{eq1}). In particular if $f(x,y)=0$ contains
an oval without any singular point of system (\ref{eq1}), this
oval is a periodic orbit of system (\ref{eq1}).
\newline

As well as invariant curves, the other objects taken into
consideration in this paper are limit cycles. A {\em limit cycle}
of system (\ref{eq1}) is an isolated periodic orbit. Let $\gamma$
be a limit cycle for system (\ref{eq1}). We say that $\gamma$ is
{\em stable} if there exists a neighborhood such that all the
orbits starting in it have $\gamma$ as $\omega$--limit set. We say
that $\gamma$ is {\em unstable} if there is a neighborhood such
that all the orbits starting in it have $\gamma$ as $\alpha$-limit
set. There might be limit cycles which are neither stable nor
unstable. These limit cycles have a neighborhood such that in the
interior of the limit cycle all the orbits have $\gamma$ as
$\omega$-limit set and in the exterior of $\gamma$ all the orbits
have $\gamma$ as $\alpha$-limit set. Or the other way round: the
orbits of the interior have $\gamma$ as $\alpha$-limit set and the
orbits in the exterior have $\gamma$ as $\omega$-limit set. In
this case, we say that $\gamma$ is {\em semi-stable}. Any limit
cycle $\gamma$ of a system (\ref{eq1}) is either stable, unstable
or semi-stable as it is stated in \cite{Perko}.
\par
A classical known result, given in the book of Perko \cite{Perko},
let us distinguish the hyperbolicity of a limit cycle. If we
consider $\gamma(t)$ a periodic orbit of system {\rm (\ref{eq1})}
of period $T$, we may compute the finite value given by the
following integral $\int_{0}^T \div(\gamma(t)) \, dt $, where
$\div(x,y)=(\partial P/\partial x) + (\partial Q/\partial y)$ is
called the {\em divergence} of system (\ref{eq1}). It can be shown
that if $\int_{0}^T \div(\gamma(t)) \, dt <0$, then $\gamma$ is a
stable limit cycle, if $\int_{0}^T \div(\gamma(t)) \, dt >0$, then
$\gamma$ is a unstable limit cycle and if $\int_{0}^T
\div(\gamma(t)) \, dt =0$, then $\gamma$ may be a stable, unstable
or semi-stable limit cycle or it may belong to a continuous band
of cycles. When the quantity $\int_{0}^T \div(\gamma(t)) \, dt$ is
different from zero, we say that the limit cycle $\gamma$ is {\em
hyperbolic}. We notice that if $\int_{0}^T \div(\gamma(t)) \, dt
\neq 0$, then the periodic orbit $\gamma$ is a limit cycle (either
stable or unstable). We are going to use this property to ensure
that a periodic orbit is a limit cycle, that is, that it does not
belong to a continuous band of cycles.
\newline

We relate limit cycles with invariants in the following way. We
assume that we have a periodic orbit $\gamma$ of system
(\ref{eq1}) which is given in an implicit way, that is, there
exists an invariant curve $f(x,y)=0$ such that $\gamma \subseteq
\{ (x,y) \, | \, f(x,y) = 0\}$. In order to have a smooth curve
$f(x,y)=0$ defining the periodic orbit, we will assume that
$\nabla f(p)\neq 0$ for any $p \in \gamma$, that is, the gradient
vector of $f(x,y)$ is different from zero in all the points of
$\gamma$. Then we have the following result stated and proved in
\cite{hyper}.
\begin{theorem}
Let us consider a system {\rm (\ref{eq1})} and $\gamma(t)$ a
periodic orbit of period $T>0$. Assume that $f:\mathcal{U}
\subseteq \mathbb{R}^2 \to \mathbb{R}$ is an invariant curve with
$\gamma \subseteq \{ (x,y) \, | \, f(x,y) = 0\}$ and let $k(x,y)$
be the cofactor of $f(x,y)$ as given in {\rm (\ref{eqci})}. We
assume that $\nabla f(p)\neq 0$ for any $p \in \gamma$. Then,
\begin{equation}
\int_{0}^T  k (\gamma(t)) \, dt \  = \  \int_{0}^T \div(\gamma(t))
\, dt. \label{eq4}
\end{equation}
\label{th4}
\end{theorem}
Hence, we have an alternative way to compute the value $\int_{0}^T
\div(\gamma(t)) \, dt$. In the family of planar polynomial
differential systems which we are considering, that is, the one
described in Theorem \ref{th1}, we will not be able to directly
compute the value $\int_{0}^T \div(\gamma(t)) \, dt$. This is due
to the fact that we are not considering a fixed system with a
concrete periodic orbit, but a family of systems each one with a
different periodic orbit and, thus, the expression of the
integrand is too general to be manipulated. By Theorem \ref{th4},
we can also compute the value $\int_{0}^T  k (\gamma(t)) \, dt$
but this integral is as much difficult as the previous one. That's
why we are going to use the fact that for any $w \in \mathbb{R}$:
\[ \int_{0}^T \div(\gamma(t)) \, dt = \int_{0}^T \div(\gamma(t)) \,
dt + w \left( \int_{0}^T \div(\gamma(t)) \, dt - \int_{0}^T
k(\gamma(t)) \, dt \right). \]  The integrand in the right hand
side of this equality will be chosen strictly positive or negative
in all the interval of integration for a suitable value of $w$.
Therefore, the value of the integral will be different from zero.
Using these steps, we will be able to prove that the oval of the
invariant curve described in Theorem \ref{th1} is a limit cycle of
the corresponding system.

\section{Proof of Theorem \ref{th1} \label{sect3}}

In order to prove this theorem, we first show that $f(x,y)=0$,
where $f(x,y):=(y-p(x) q(x))^2-p(x)$, is an invariant algebraic
curve of system (\ref{eq2}). Straightforward computations show
that:
\[ \begin{array}{l}
\displaystyle y \left( \frac{\partial f}{\partial x} \right) +
\left[ \left(\frac{3}{2} \, q(x) p'(x) + p(x) q'(x) \right) y -
\displaystyle \frac{p'(x)}{2} \left( p(x) q(x)^2-1 \right) \right]
\left( \frac{\partial f}{\partial y} \right) = \vspace{0.2cm} \\
\displaystyle = q(x) p'(x) f(x,y), \end{array} \] and, thus, we
have that $f(x,y)=0$ is an invariant algebraic curve for system
(\ref{eq2}) with cofactor $k(x,y):= q(x) p'(x)$.
\par
Since $p(x_e)=p(x_d)=0$ for the values $x_e<0$ and $x_d>0$ and
$p(x) >0$ for $x \in (x_e,x_d)$, we deduce that $f(x,y)=0$ has an
oval in the band $x_e \leq x \leq x_d$ surrounding the origin of
coordinates, which can be parameterized in two parts by:
\begin{equation}
x(\tau) = \tau, \quad y_{\pm}(\tau) = p(\tau) q(\tau) \pm
\sqrt{p(\tau)}, \label{eq5} \end{equation} with $\tau \in
(x_e,x_d)$. We are going to prove that this oval does not contain
any singular point of system (\ref{eq2}), and then, we will have
that it defines a periodic orbit of the system. The singular
points of system (\ref{eq2}) have coordinates of the form $(a,0)$
where the value $a$ is a root of the polynomial $p'(x) (p(x)
q(x)^2-1)$. We have that $f(a,0)= p(a) (p(a) q(a)^2-1)$ and $(p(a)
q(a)^2-1)$ is different from zero in all the closed interval $a
\in [x_e, x_d]$ by the hypothesis that $p(x) q(x)^2 \neq 1$ for
all $x \in (x_e,x_d)$ and $p(x_e)=p(x_d)=0$.
\par
Here we notice that the assumption $p(x) q(x)^2 \neq 1$ for $x \in
(x_e,x_d)$ is necessary for the oval of $f(x,y)=0$ to be a limit
cycle. In \cite{Abdelkader}, this assumption is not given. We
notice that an oval of an invariant algebraic curve of a system
may contain singular points of the system, and in such a case, it
is not even a periodic orbit.
\par
Since $x_e$ and $x_d$ are simple zeroes of $p(x)$ and $p(x)>0$ for
all $x$ in the interval $(x_e,x_d)$ we have that there is no
singular point of system (\ref{eq2}) on the oval given by
$f(x,y)=0$ and parameterized by (\ref{eq5}). From this fact and
that $f(x,y)=0$ is an invariant algebraic curve of the system, we
deduce that this oval is a periodic orbit of system (\ref{eq1}).
We denote this periodic orbit by $\gamma$ for the rest of the
proof. We note that we do not know the parameterization of
$\gamma$ as explicit solution of system (\ref{eq2}), that is, we
do not know the periodic function $\gamma(t):=(\gamma_1(t),
\gamma_2(t))$ such that $d\gamma_1(t)/dt = \gamma_2(t)$ and
\begin{eqnarray*} \frac{d \gamma_2(t)}{dt} & = & \displaystyle
\left(\frac{3}{2} \, q(\gamma_1(t)) p'(\gamma_1(t)) +
p(\gamma_1(t)) q'(\gamma_1(t)) \right) \gamma_2(t)  \\ & & -
\displaystyle \frac{p'(\gamma_1(t))}{2} \left( p(\gamma_1(t))
q(\gamma_1(t))^2-1 \right), \end{eqnarray*} for all $t \in
\mathbb{R}$. We do neither know its period $T>0$ but we have been
able to show its existence by using the invariant algebraic curve
$f(x,y)=0$ and its properties in relation with system (\ref{eq2}).
\par
Finally, we need to prove that the periodic orbit $\gamma$ is a
hyperbolic limit cycle. To do so, we are going to show that the
value of the integral $\int_{0}^{T} \div(\gamma(t)) \, dt$ is
different from zero. Since we do not know $\gamma(t)$ nor the
period $T$, we use the parameterization of the oval $\gamma$ given
in (\ref{eq5}). In order to get the correct sign of the integral
$\int_{0}^{T} \div(\gamma(t)) \, dt$, we need to know the sense of
the flow over $\gamma$. We take the point of coordinates
$(x_d,0)$, which belongs to $\gamma$, and we have that the vector
field defined by system (\ref{eq2}) on that point is $(0,
p'(x_d)/2)$ because $p(x_d)=0$. Since $p(x_d)=0$, $p'(x_d) \neq 0$
and $p(x)>0$ in the interval $(x_e,x_d)$, we deduce that
$p'(x_d)<0$. Hence, the sense of the flow over $\gamma$ is
clockwise. We can write the following equality, using the
parameterization (\ref{eq5}):
\begin{eqnarray*} \displaystyle \int_{0}^{T} \div(\gamma(t)) \, dt & = &
\int_{x_e}^{x_d} \frac{\div(x(\tau),y_{+}(\tau))}{y_{+}(\tau)} \,
d \tau + \int_{x_d}^{x_e}
\frac{\div(x(\tau),y_{-}(\tau))}{y_{-}(\tau)} \, d \tau \\ & = &
\displaystyle \int_{x_e}^{x_d} \left(
\frac{\div(x(\tau),y_{+}(\tau))}{y_{+}(\tau)} -
\frac{\div(x(\tau),y_{-}(\tau))}{y_{-}(\tau)} \right) \, d \tau .
\end{eqnarray*}
We note that the divergence of system (\ref{eq2}) is $\div(x,y)= 3
\, q(x) p'(x)/2 + p(x) q'(x)$, and substituting this expression in
the former equality, we get:
\[ \int_{0}^{T} \div(\gamma(t)) \, dt \, = \, \int_{x_e}^{x_d}
\frac{ - 3 q(\tau)\, p'(\tau) - 2 p(\tau) \, q'(\tau)}{(p(\tau)
q(\tau)^2-1) \sqrt{p(\tau)}}\, d \tau. \] This integral is well
defined because we are assuming that $(p(\tau) q(\tau)^2-1)
p(\tau)$ is different from zero for $\tau \in (x_e,x_d)$. However,
we are not able to distinguish if its value is positive, negative
or zero. By using the same reasonings, we can write the following
equality:
\[ \int_{0}^{T} k(\gamma(t)) \, dt \, = \, \int_{x_e}^{x_d}
\frac{ - 2 q(\tau)\, p'(\tau)}{(p(\tau) q(\tau)^2-1)
\sqrt{p(\tau)}}\, d \tau. \] Using Theorem \ref{th4}, we have
that, for any value of $w \in \mathbb{R}$:
\begin{eqnarray*}
\displaystyle \int_{0}^{T} \div(\gamma(t)) \, dt & = &
\displaystyle \int_{0}^{T} \div(\gamma(t)) \, dt + w \left(
\displaystyle \int_{0}^{T} \div(\gamma(t))\, dt - \displaystyle
\int_{0}^{T} k(\gamma(t)) \, dt \right) \\ & = & \displaystyle
\int_{x_e}^{x_d} \frac{ -(w+3) q(\tau)p'(\tau) - 2 (1+w) p(\tau)
q'(\tau)}{(p(\tau) q(\tau)^2-1) \sqrt{p(\tau)}} \, d \tau.
\end{eqnarray*} Taking $w=-3$, we get:
\begin{equation}
\displaystyle \int_{0}^{T} \div(\gamma(t)) \, dt \, = \,
\int_{x_e}^{x_d} \frac{ 4 p(\tau) q'(\tau)}{(p(\tau) q(\tau)^2-1)
\sqrt{p(\tau)}} \, d \tau \label{eq6}
\end{equation}
The hypothesis of Theorem \ref{th1} on the sign of the polynomials
$p(x)$ and $q(x)$ in the interval $x \in (x_e,x_d)$ are $p(x)>0$,
$p(x) q(x)^2 \neq 1$ and $q'(x) \neq 0$. Therefore the integrand
of the right hand side of (\ref{eq6}) is strictly positive or
negative in all the interval $\tau \in (x_e,x_d)$. We deduce that
the value of the integral cannot be zero and, hence, the periodic
orbit $\gamma$ is a hyperbolic limit cycle as we wanted to show.
\par
We would also like to characterize if this limit cycle is stable
or unstable, so we are going to study the sign of the integrand in
the right hand side of (\ref{eq6}). Since $p(x) q(x)^2 \neq 1$ for
$x \in (x_e,x_d)$ and $p(x_d) q(x_d)^2-1=-1$ (because $p(x_d)=0$),
we deduce that $p(x) q(x)^2-1 <0$ for $x \in (x_e,x_d)$.
Therefore, we have that the integrand in the right hand side of
(\ref{eq6}) is strictly positive if $q'(\tau) < 0$ for all $\tau
\in (x_e,x_d)$ and strictly negative if $q'(\tau) > 0$ for all
$\tau \in (x_e,x_d)$. We can state that the hyperbolic limit cycle
$\gamma$ is stable if $q'(0)>0$ and unstable if $q'(0)<0$. \bbox
\newline

We also note that in the work \cite{Abdelkader}, the expression of
an example of a more general limit cycle for a family of planar
polynomial differential systems is pretended to be given. In fact,
we are going to show that in the case that the oval of this
example is a limit cycle, we are in the same family of systems as
written in (\ref{eq2}), that is, the one described in Theorem
\ref{th1}.
\par
In \cite{Abdelkader} the following planar polynomial differential
system is given as an example of a more general family of systems
with an explicit limit cycle.
\begin{equation}
\begin{array}{lll}
\displaystyle \dot{x} & = & \displaystyle y, \vspace{0.2cm} \\
\displaystyle \dot{y} & = & \displaystyle \left\{ p'(x) \left[
(m+r) h(x) p(x)^{r-1} + (m+1) q(x) \right] + h'(x) p(x)^r \right.
\vspace{0.2cm} \\
\displaystyle & & \left. \displaystyle + p(x) q'(x) \right\} y  -
\displaystyle m p(x) p'(x) \left( \left[ h(x) p(x)^{r-1} + q(x)
\right]^2 - p(x)^{2m-2} \right),
\end{array} \label{eq7}
\end{equation}
where $p(x)$ is an even polynomial, $q(x)$ and $h(x)$ are odd
polynomials, $m=n+1/2$, $n$ is an integer number with $n \geq 0$
and $r$ is an integer number with $r \geq 2$. Moreover, it is
assumed that $p(0)>0$ and there exists a value $X>0$ such that
$p(X)=0$, $p(x)<0$ for all $x>X$ and $q(0)=0$. \par Some
straightforward computations show that system (\ref{eq7}) exhibits
the invariant algebraic curve $f(x,y)=0$ with $f(x,y):=
(y-h(x)p(x)^r - q(x) p(x))^2-p(x)^{2n+1}$ and with cofactor
$k(x,y):= (2n+1) p'(x) \left[ h(x) p(x)^{r-1} + q(x) \right]$. We
have that $f(x,y)=0$ has an oval in the band $-X \leq x \leq X$
which can be parameterized by:
\[ x(\tau)=\tau, \ y(\tau)= h(\tau) p(\tau)^r + p(\tau) q(\tau)
\pm p(\tau)^n \sqrt{p(\tau)}. \] In order to show that this oval
is a periodic orbit of system (\ref{eq7}), we only need to show
that it does not contain any singular point of the system. The
singular points of system (\ref{eq7}) in the band $|x|\leq X$ are
of the form $(a,0)$ where $a$ is a root of the polynomial $p(x)
p'(x) \left( \left[ h(x) p(x)^{r-1} + q(x) \right]^2 - p(x)^{2n-1}
\right)$, because $m=n+1/2$. We notice that, unless $n=0$, the
points with coordinates $(-X,0)$ and $(X,0)$ are singular points
of the system which are contained in the oval of $f(x,y)=0$.
Therefore, if $n>0$, we have that the oval of $f(x,y)=0$ cannot be
a limit cycle. If $n=0$, we can consider the polynomial
$\tilde{q}(x):= q(x) + h(x) p(x)^{r-1}$ and we have that system
(\ref{eq7}) coincides with system (\ref{eq2}) with polynomials
$p(x)$ and $\tilde{q}(x)$. Therefore, this is not an example of a
more general limit cycle.


\begin{thebibliography}{99}

\bibitem{Abdelkader} {\sc M.A. Abdelkader},{\it \ Relaxation oscillators
with exact limit cycles.} J. Math. Anal. Appl. {\bf 218} (1998)
308--312.

\bibitem{hyper} {\sc H. Giacomini and M. Grau}, {\it \ On the stability of
limit cycles for planar differential systems}, J. Differential
Equations {\bf 213} (2005), 368--388.

\bibitem{Odani} {\sc K. Odani}, {\it \ The limit cycle of the van der Pol
equation is not algebraic}, J. Differential Equations {\bf 115}
(1995), 146--152.

\bibitem{Perko} {\sc L. Perko},{\it \ Differential equations
and dynamical systems.} Third edition. Texts in Applied
Mathematics, {\bf 7}. Springer-Verlag, New York, 2001.

\end{thebibliography}
\end{document}